\documentclass{amsart}
\usepackage{amsmath,amssymb}
\usepackage{yhmath}

\newcommand{\bdis}{\begin{displaymath}}
\newcommand{\edis}{\end{displaymath}}
\newcommand{\be}{\begin{equation}}
\newcommand{\ee}{\end{equation}}
\newcommand{\mbb}{\mathbb}
\newcommand{\mcal}{\mathcal}

\newcommand{\vp}{\varphi}

\newcommand{\zf}{\zeta\left(\frac{1}{2}+it\right)}

\DeclareMathOperator{\re}{Re}

\theoremstyle{definition}

\theoremstyle{remark}
\newtheorem{remark}[]{Remark}

\newtheorem*{mydef11}{{\bf Theorem 1}}

\newtheorem*{mydef12}{{\bf Theorem 2}}

\newtheorem*{mydef13}{{\bf Theorem 3}}

\newtheorem*{mydef41}{{\bf Corollary 1}}

\newtheorem*{mydef42}{{\bf Corollary 2}}

\newtheorem*{mydef43}{{\bf Corollary 3}}

\newtheorem*{mydef44}{{\bf Corollary 4}}

\newtheorem*{mydef45}{{\bf Corollary 5}}

\newtheorem*{mydef91}{{\bf Formula1}}

\newtheorem*{mydef92}{{\bf Formula2}}

\newtheorem*{mydef93}{{\bf Formula3}}

\numberwithin{equation}{section}

\begin{document}

\title{Jacob's ladders, existence of almost linear increments of the Hardy-Littlewood integral and new types of multiplicative laws}{}

\author{Jan Moser}

\address{Department of Mathematical Analysis and Numerical Mathematics, Comenius University, Mlynska Dolina M105, 842 48 Bratislava, SLOVAKIA}

\email{jan.mozer@fmph.uniba.sk}


\begin{abstract}
In this paper we prove that there is a continuum set of increments with some minimal structure for the Hardy - Littlewood integral. The result implies a number of new properties of the Hardy - Littlewood integral. 
\end{abstract}
\maketitle

\section{Introduction} 

\subsection{}    

In this paper we shall prove some new results of completely new type for certain class of increments of the classical Hardy - Littlewood integral\footnote{See \cite{3}.}
\be \label{1.1} 
J(T)=\int_0^T\left|\zf\right|^2{\rm d}t. 
\ee 
We use for this purpose the theory of Jacob's ladders introduced in our papers \cite{7} -- \cite{10}. Namely it is true: for every sufficiently big $T>0$ and for every fixed $k\in\mbb{N}$ there is a sequence 
\be \label{1.2} 
\{ \overset{r}{T}(T)\}_{r=1}^k:\ T<\overset{1}{T}(T)<\overset{2}{T}(T)<\dots<\overset{k}{T}(T)
\ee 
that the following formula 
\be \label{1.3} 
\int_{\overset{r-1}{T}}^{\overset{r}{T}}\left|\zf\right|^2{\rm d}t\sim (1-c)\overset{r-1}{T},\ T\to\infty 
\ee 
holds true, where $c$ is the Euler's constant. 

Since (see (\ref{1.1})) 
\be \label{1.4} 
\int_{\overset{r-1}{T}}^{\overset{r}{T}}\left|\zf\right|^2{\rm d}t=J(\overset{r}{T})-J(\overset{r-1}{T})=\Delta J(\overset{r-1}{T},\overset{r}{T}),\ r=1,\dots,k
\ee 
then the set of integrals in (\ref{1.3}) represents the set of correspondings increments of the Hardy-Littlewood integral (\ref{1.1}). 

\begin{remark}
Consequently, we have constructed vector function of increments of the Hardy - Littlewood integral 
\be \label{1.5} 
\int_{\overset{r-1}{T}}^{\overset{r}{T}}\left|\zf\right|^2{\rm d}t,\ r=1,\dots,k,\ T\in (T_0,\infty), 
\ee  
with positive and sufficiently big $T_0$ that fulfils the following property: 
\be \label{1.6} 
\lim_{T\to\infty}\frac{1}{\overset{r-1}{T}}\int_{\overset{r-1}{T}}^{\overset{r}{T}}\left|\zf\right|^2{\rm d}t=1-c,\ r=1,\dots,k. 
\ee 
\end{remark} 

\subsection{} 

Now, we shall present some strange formulas in the theory of the Riemann zeta-functions following from (\ref{1.3}). 

\begin{itemize}
	\item[(A)] Multiplicative formulas of the type: 
	\be \label{1.7} 
	\begin{split} 
	& \int_{T\times\overset{1}{T}\times \cdots\times\overset{k-1}{T}}^{[T\times\overset{1}{T}\times \cdots\times\overset{k-1}{T}]^1}\left|\zf\right|^2{\rm d}t\sim  \\ 
	& 
	\frac{1}{(1-c)^{k-1}}\prod_{r=1}^k\int_{\overset{r-1}{T}}^{\overset{r}{T}}\left|\zf\right|^2{\rm d}t,\ T\to\infty, 
	\end{split} 
	\ee 
	where we use the following notation 
	\be \label{1.8}  
	[T\times\overset{1}{T}\times \cdots\times\overset{k-1}{T}]^1=\overset{1}{\wideparen{T\times\overset{1}{T}\times \cdots\times\overset{k-1}{T}}}. 
	\ee  
	\item[(B)]Additive formula of the type: 
	\be \label{1.9} 
	\int_{T+\overset{1}{T}+ \cdots+\overset{k-1}{T}}^{[T+\overset{1}{T}+ \cdots+\overset{k-1}{T}]^1}\left|\zf\right|^2{\rm d}t\sim 
	\sum_{r=1}^k \int_{\overset{r-1}{T}}^{\overset{r}{T}}\left|\zf\right|^2{\rm d}t,\ T\to\infty. 
	\ee 
	\item[(C)] And mixed formulas, for example: 
	\be \label{1.10} 
	\begin{split}
	& \int_{T+\overset{1}{T}\overset{2}{T}}^{[T+\overset{1}{T}\overset{2}{T}]^1}\left|\zf\right|^2{\rm d}t\sim \\ 
	& \int_{T}^{\overset{1}{T}}\left|\zf\right|^2{\rm d}t+ \\ 
	& \frac{1}{1-c}\int_{\overset{1}{T}}^{[\overset{1}{T}]^1}\left|\zf\right|^2{\rm d}t\times 
	\int_{\overset{2}{T}}^{[\overset{2}{T}]^1}\left|\zf\right|^2{\rm d}t,\ T\to\infty. 
	\end{split}
	\ee 
\end{itemize} 

\begin{remark}
The formula (\ref{1.7}) can be rewritten in the form 
\be \label{1.11} 
\begin{split}
& \int\int_{\square_k}\dots\int\prod_{r=1}^k \left|\zeta\left(\frac 12+i(x_{r-1}+\overset{r-1}{T})\right)\right|^2{\rm d}x_0{\rm d}x_1\dots{\rm d}x_{k-1}\sim \\ 
& \int_{T\times\overset{1}{T}\times \cdots\times\overset{k-1}{T}}^{[T\times\overset{1}{T}\times \cdots\times\overset{k-1}{T}]^1}\left|\zf\right|^2{\rm d}t,\ T\to\infty, 
\end{split}
\ee  
where 
\be \label{1.12} 
\square_k:\ 0\leq x_{r-1}\leq \overset{r}{T}-\overset{r-1}{T},\ r=1,\dots,k,\ \square_k\subset\mbb{R}^k. 
\ee 
\end{remark} 

\begin{remark}
Our formula (\ref{1.7}) establishes a nonlinear asymptotic dependence of each of the integrals in the set 
\be \label{1.13} 
\left\{
\int_{T\times\overset{1}{T}\times \cdots\times\overset{k-1}{T}}^{[T\times\overset{1}{T}\times \cdots\times\overset{k-1}{T}]^1}\left|\zf\right|^2{\rm d}t,\ 
\int_{T}^{\overset{1}{T}}\left|\zf\right|^2{\rm d}t,\dots, 
\int_{\overset{k-1}{T}}^{\overset{k}{T}}\left|\zf\right|^2{\rm d}t
\right\}
\ee 
on remaining elements. On the other hand, our formula (\ref{1.9}) establishes a linear dependence between the integrals in the set 
\be \label{1.14} 
\left\{
\int_{T+\overset{1}{T}+ \cdots+\overset{k-1}{T}}^{[T+\overset{1}{T}+ \cdots+\overset{k-1}{T}]^1}\left|\zf\right|^2{\rm d}t, \ 
\int_{T}^{\overset{1}{T}}\left|\zf\right|^2{\rm d}t, \dots,
\int_{\overset{k-1}{T}}^{\overset{k}{T}}\left|\zf\right|^2{\rm d}t
\right\}. 
\ee 
\end{remark} 

\begin{remark} 
Main distances of the segments of integration in (\ref{1.7}) are given by: 
\be \label{1.15} 
\begin{split}
& \rho_k\{[T\times\overset{1}{T}\times \cdots\times\overset{k-1}{T},[T\times\overset{1}{T}\times \cdots\times\overset{k-1}{T}]^1];[\overset{k-1}{T},\overset{k}{T}]\}= \\ 
& T\times\overset{1}{T}\times \cdots\times\overset{k-1}{T}-\overset{k}{T}>T^k-(1+\epsilon)T>0.9\times T^k. 
\end{split}
\ee 
Comp. (3.3) with $0<\epsilon$ being sufficiently small. Therefore they are big with big $k$. For example 
\be \label{1.16} 
\rho_{257}>0.9\times T^{257}. 
\ee 
\end{remark} 

\subsection{} 

Next, we have obtained also a nonlinear formula that generates corresponding Hardy - Littlewood integral by means of the increments, (\ref{1.3}) and (\ref{1.4}). Namely the following equality holds true: 
\be \label{1.17} 
\begin{split}
& \int_0^{\overset{r-1}{T}}\left|\zf\right|^2{\rm d}t\sim \\ 
& \frac{1}{1-c}\int_{\overset{r-1}{T}}^{\overset{r}{T}}\left|\zf\right|^2{\rm d}t\times 
\ln\left\{ 
\frac{e^{\Lambda_1}}{1-c}\int_{\overset{r-1}{T}}^{\overset{r}{T}}\left|\zf\right|^2{\rm d}t
\right\}, \\ 
& r=1,\dots,k,\ T\to\infty, 
\end{split}
\ee  
where 
\be \label{1.18} 
\Lambda_1=2c-1-\ln 2\pi. 
\ee  

\begin{remark}
Our formula (\ref{1.17}) expresses the following property: \emph{smaller} integral generates \emph{bigger} one. 
\end{remark} 

\begin{remark}
As a consequence, our result (\ref{1.3}) and its implications labelled (\ref{1.7}) -- (\ref{1.11}) and (\ref{1.17})\footnote{And also some others, we shall present later.} express new kind of phenomena connected with the Hardy - Littlewood integral (\ref{1.1}) on distances 
\bdis 
(1-c)\pi(T), 
\edis  
where, as usually, the symbol $\pi(T)$ stands for the prime-counting function. 
\end{remark} 

\section{Brief history of the Hardy - Littlewood integral (1918)} 

\subsection{} 

The origin of the Hardy - Littlewood integral lies in the asymptotic formula 
\be \label{2.1} 
\int_0^T\left|\zf\right|^2{\rm d}t\sim T\ln T,\ T\to\infty, 
\ee  
see \cite{3}, pp. 122, 151 -- 156. Littlewood formula \cite{5} ($0<\delta$ being sufficiently small): 
\be \label{2.2} 
\int_0^T\left|\zf\right|^2{\rm d}t=T\ln T-(1+\ln 2\pi)T+\mcal{O}(T^{\frac 34+\delta}) 
\ee 
was the next step. And finally, Ingham proved the following formula \cite{6}: 
\be \label{2.3} 
\int_0^T\left|\zf\right|^2{\rm d}t=T\ln T-(1+\ln 2\pi - 2c)T+\mcal{O}(T^{\frac 12}\ln T), 
\ee  
that contains complete set of main terms. 

\begin{remark}
We see that the first Ingham's contribution to the formula (\ref{2.3}) is in addition of the term $-2cT$ that was missing in the previous Littlewood formula (\ref{2.2}). The second contribution of Ingham is quite essential improvement of the error term in Littlewood's formula, namely by $33.3\%$. 
\end{remark} 

\begin{remark}
Let us remind also some of other considerable refinements 
\be \label{2.4} 
\mcal{O}(T^{\frac{5}{12}+\delta}),\ 16.6\%;\ \mcal{O}(T^{\frac 13+\delta}),\ 20\% 
\ee 
of the error term in (\ref{2.3}) due to Titschmarsh \cite{11} and Balasubramanian \cite{1}, correspondingly.  
\end{remark} 

We can include previous results into the following statement. 

\begin{mydef91}
\be\label{2.5} 
\int_0^T\left|\zf\right|^2{\rm d}t=T\ln T-(1+\ln 2\pi - 2c)T+R(T). 
\ee 
We will refer to (\ref{2.5}) as to Hardy-Littlewood-Ingham formula with the error term $R(T)$. 
\end{mydef91} 

\subsection{} 

Further, let us remind the fundamental Good's $\Omega$-theorem, \cite{2}: 
\be \label{2.6}
R(T)=\Omega(T^{\frac 14}),\ T\to\infty. 
\ee  

\begin{remark}
It is possible to write down 
\be \label{2.7} 
R(T)=\mcal{O}(T^{a+\delta}),\ \frac 14\leq a\leq \frac 13, 
\ee 
for the error term. In this form $a$ denotes the least valid exponent for which the estimate (\ref{2.7}) holds true and $a=\frac 14$ represents the Good's $\Omega$-theorem and $a=\frac 13$ represents the main Balasubramanian's exponent. 
\end{remark} 

\begin{remark}
It is fully sufficient to use the exponent $\frac 13$ despite there is also the second Balasubramanian's exponent that is equal to $\frac{27}{82}$, that is an improvement by $1.18\%$, see \cite{1}. 
\end{remark} 

\section{Main result as the next step in the evolution of the Hardy - Littlewood formula (1918)} 

\subsection{} 

We use the following notions from our papers \cite{7} -- \cite{10}: 
\begin{itemize}
	\item[(a)] The Jacob's ladder $\vp_1(T)$, 
	\item[(b)] direct iterations of Jacob's ladder: 
	\bdis 
	\begin{split}
		& \vp_1^0(T)=T,\ \vp_1^1(T)=\vp_1(T),\ \vp_1^2(T)=\vp_1(\vp_1(T)),\ \dots, \\ 
		& \vp_1^r(T)=\vp_1(\vp_1^{r-1}(T)),\ r=1,\dots,k, 
	\end{split}
	\edis 
	\item[(c)] reverse iterations of Jacob's ladder (by means of the inverse function $\vp_1^{-1}(T)$): 
	 \be \label{3.1}  
	 \begin{split}
	 	& \vp_1^{-1}(T)=\overset{1}{T},\ \vp_1^{-2}(T)=\vp_1^{-1}(\overset{1}{T})=\overset{2}{T},\dots , \\ 
	 	& \vp_1^{-r}(T)=\vp_1^{-1}(\overset{r-1}{T})=\overset{r}{T},\ r=1,\dots,k, 
	 \end{split}
	\ee 
	that is, for example, 
	\be \label{3.2} 
	\vp_1(\overset{r}{T})=\overset{r-1}{T}
	\ee 
	for every fixed $k\in\mbb{N}$ and every sufficiently big $T>0$. 
\end{itemize} 

Next, we use also the following properties:\footnote{See \cite{9}} 
\be \label{3.3} 
\begin{split}
& \overset{r}{T}-\overset{r-1}{T}\sim (1-c)\pi(\overset{r}{T}),\ \pi(\overset{r}{T})\sim \frac{\overset{r}{T}}{\ln \overset{r}{T}},\ r=1,\dots,k, \\ 
& \overset{0}{T}=T<\overset{1}{T}(T)<\overset{2}{T}(T)<\dots<\overset{k}{T}(T), \\ 
& T\sim\overset{1}{T}\sim\overset{2}{T}\sim\dots\sim\overset{k}{T},\ T\to\infty. 
\end{split}
\ee 

\begin{remark}
The asymptotic behaviour of the points 
\bdis 
\{ T,\overset{1}{T},\dots,\overset{k}{T}\}
\edis 
is as follows: As $T\to\infty$ these points recede unboundedly each from other and all together are receding to infinity. Hence, the set of these points behaves at $T\to\infty$ as one-dimensional Friedmann-Hubble expanding Universe. 
\end{remark} 

\subsection{} 

Here we formulate our main result. 

\begin{mydef11} 
For every fixed $k\in\mbb{N}$ and every sufficiently big $T>0$ we have: 
\be \label{3.4} 
\begin{split}
& \int_{\overset{r-1}{T}}^{\overset{r}{T}}\left|\zf\right|^2{\rm d}t=(1-c)\overset{r-1}{T}+\mcal{O}(T^{a+\delta}), \\ 
& 1-c\approx 0.42,\ r=1,\dots,k,\ T\to\infty 
\end{split}
\ee  
where, by (\ref{3.1}): 
\be \label{3.5} 
\overset{r}{T}=\overset{r}{T}(T)=\vp_1^{-r}(T). 
\ee 
\end{mydef11}  

\begin{remark}
The existence of linear main terms in (\ref{3.4}), linear with respect to variables (\ref{3.2}), is a new phenomenon in the theory of the Riemann's $\zf$-function. For example, it is true\footnote{Remark 9, $a=\frac 13$.} that 
\be \label{3.6} 
\int_{\overset{r-1}{T}}^{\overset{r}{T}}\left|\zf\right|^2{\rm d}t=(1-c)\overset{r-1}{T}+\mcal{O}(T^{\frac 13+\delta}),\ T\to\infty. 
\ee 
\end{remark} 

\begin{remark} 
Let us notice explicitly, that our formula (\ref{3.4}) inherits every true estimate for its error term from HLI formula (\ref{2.5}).\footnote{Comp. (\ref{2.7}) and Remark 8.} 
\end{remark}  

\section{The first class of consequences of (\ref{3.4})} 

\subsection{} 

At the first place the following statement is obtained by (\ref{3.4}). 

\begin{mydef41}
\be \label{4.1} 
\begin{split}
& \frac{1}{\overset{r-1}{T}}\int_{\overset{r-1}{T}}^{\overset{r}{T}}\left|\zf\right|^2{\rm d}t=1-c+\mcal{O}(T^{-1+a+\delta}),\\ 
& r=1,\dots,k,\ T\to\infty, 
\end{split}
\ee  
where (see Remark 9) 
\be \label{4.2} 
\mcal{O}(T^{-1+a+\delta})=\mcal{O}(T^{-\frac 23+\delta}). 
\ee 
\end{mydef41} 

\begin{remark}
The formula (\ref{4.1}) expresses certain asymptotic conservation law controlling chaotic oscillations of the Riemann's $\zf$-function. 
\end{remark} 

Let us remind in connection with chaotic oscillations of the Riemann's zeta-function on the critical line the following: 

\begin{itemize}
	\item[(A)] Spectral form of the Riemann-Siegel formula\footnote{Comp. \cite{10}, subsection 1.4.} 
	\be \label{4.3} 
	\begin{split}
	& Z(t)=\sum_{n\leq \tau(x)}\frac{2}{\sqrt{n}}\cos\{t\omega_n(x)+\psi(x)\}+Q(x), \\ 
	& \tau(x)=\sqrt{\frac{x}{2\pi}},\ Q(x)=\mcal{O}(x^{-\frac 14}),\ t\in [x,x+V], \\ 
	& V=(0,x^{\frac 14}], x\in [T,\overset{k}{T}], \\ 
	& |Z(t)|=\left|\zf\right|, 
	\end{split}
	\ee 
	where the functions 
	\bdis 
	\frac{2}{\sqrt{n}}\cos\{t\omega_n(x)+\psi(x)\}
	\edis  
	are Riemann's oscillators with: 
	\begin{itemize}
		\item[(a)] Amplitude 
		\bdis 
		\frac{2}{\sqrt{n}}, 
		\edis  
		\item[(b)] Incoherent local phase constant 
		\bdis 
		\psi(x)=-\frac x2-\frac{\pi}{8}, 
		\edis 
		\item[(c)] Nonsynchronized local times 
		\bdis 
		t\in [x,x+V], 
		\edis  
		\item[(d)] Local spectrum of cyclic frequencies 
		\bdis 
		\{\omega_n(x)\}_{n\leq \tau(x)},\ \omega_n(x)=\ln\frac{\tau(x)}{n}. 
		\edis 
	\end{itemize} 
	\item[(B)] The graph of the function $|\zf|$ after transformation 
	\bdis 
	Z(t)\to|Z(t)|=\left|\zf\right|, 
	\edis  
	in the neighbourhood of the first Lehmer's pair of zeroes of $Z(t)$, see \cite{4}, pp. 296, 297, for example. 
\end{itemize} 

\subsection{} 

Next, let us remind the Euler's $c$-constant definition 
\bdis 
c=\lim_{n\to\infty}\left\{ 1+\frac 12+\dots+\frac 1n - \ln n\right\} \approx 0.5772157 
\edis 
and, for example, Dirichlet's formula for the same constant 
\be \label{4.4} 
c=\int_0^\infty \left\{ \frac{1}{1+t}-e^{-t}\right\}{\rm d}t. 
\ee  

Now we obtain the following result in this direction. 

\begin{mydef42}
\be \label{4.5} 
c=1-\lim_{T\to\infty}\frac{1}{\overset{r-1}{T}}\int_{\overset{r-1}{T}}^{\overset{r}{T}}\left|\zf\right|^2{\rm d}t, 
\ee 
or 
\be \label{4.6} 
\bar{c}=\lim_{T\to\infty}\frac{1}{\overset{r-1}{T}}\int_{\overset{r-1}{T}}^{\overset{r}{T}}\left|\zf\right|^2{\rm d}t
\ee 
for \emph{the complementary} Euler's constant $\bar{c}=1-c$. 
\end{mydef42} 

\begin{remark}
Our formula (\ref{4.5}) gives new expression for the Euler's $c$ and this expression is generated by the Riemann's $\zf$-function together with the reverse iterations of the Jacob's ladder, see (\ref{3.1}). 
\end{remark} 

\subsection{} 

Since 
\be \label{3.7} 
\frac{1}{\overset{r}{T}-\overset{r-1}{T}}=\mcal{O}\left(\frac{\ln T}{T}\right), 
\ee  
see (\ref{3.3}), then we have the following. 

\begin{mydef43}
\be \label{4.8} 
\begin{split}
& \frac{1}{\overset{r}{T}-\overset{r-1}{T}}\int_{\overset{r-1}{T}}^{\overset{r}{T}}\left|\zf\right|^2{\rm d}t= \\ 
& (1-c)\frac{\overset{r-1}{T}}{\overset{r}{T}-\overset{r-1}{T}}+\mcal{O}(T^{-1+a+\delta}\ln T),\ r=1,\dots,k,\ T\to\infty. 
\end{split}
\ee 
\end{mydef43} 

\begin{remark}
Another new phenomenon is stated here. Namely existence of rational expressions for the main members of the mean values in variables (\ref{3.2}). 
\end{remark} 

\subsection{} 

Since, see (\ref{3.3})  
\be \label{4.9} 
\overset{r}{T}-\overset{r-1}{T}\geq (1-o(1))(1-c)\frac{T}{\ln T},\ T\to\infty, 
\ee 
and (see (\ref{3.4}) and (\ref{3.6})) 
\be \label{4.10} 
T^{a+\delta}<T^{\frac 13+\delta}, 
\ee  
then we have the following asymptotic formula. 

\begin{mydef44}
\be \label{4.11} 
\begin{split}
& \int_{\overset{r}{T}}^{\overset{r+1}{T}}\left|\zf\right|^2{\rm d}t-\int_{\overset{r-1}{T}}^{\overset{r}{T}}\left|\zf\right|^2{\rm d}t= \\ 
& (1-c)(\overset{r}{T}-\overset{r-1}{T})+\mcal{O}(T^{a+\delta}),\ r=1,\dots,k,\ T\to\infty. 
\end{split}
\ee 
\end{mydef44} 

\section{The second class of consequences of (\ref{3.4})} 

\subsection{} 

Since (see (\ref{3.3}) and (\ref{3.6})) 
\be \label{5.1} 
\int_{\overset{r-1}{T}}^{\overset{r}{T}}\left|\zf\right|^2{\rm d}t=(1-c)\overset{r-1}{T}\left\{ 1+\mcal{O}(T^{-\frac 23+\delta})\right\}, 
\ee  
then 
\be \label{5.2} 
\prod_{r=1}^k\int_{\overset{r-1}{T}}^{\overset{r}{T}}\left|\zf\right|^2{\rm d}t=\left\{ 1+\mcal{O}(T^{-\frac 23+\delta}\right\}\times 
\prod_{r=1}^{k}(1-c)\overset{r-1}{T}.
\ee 

Next, we have by (\ref{5.1}) with $r=1$ 
\be \label{5.3} 
\int_T^{\overset{1}{T}}\left|\zf\right|^2{\rm d}t=(1-c)T\left\{ 1+\mcal{O}(T^{-\frac 23+\delta})\right\}
\ee 
and we put herein the following\footnote{See (\ref{3.3}).} 
\be \label{5.4} 
\begin{split}
& T\to T\times \overset{1}{T}\times \dots\times \overset{k-1}{T}\sim (T)^k, \\ 
& \overset{1}{T}\to [T\times \overset{1}{T}\times \dots\times \overset{k-1}{T}]^1, 
\end{split}
\ee  
with the result 
\be \label{5.5} 
\begin{split}
& \int_{T\times \overset{1}{T}\times \dots\times \overset{k-1}{T}}^{[T\times \overset{1}{T}\times \dots\times \overset{k-1}{T}]^1}\left|\zf\right|^2{\rm d}t= \\ 
& \left\{ 1+\mcal{O}(T^{-\frac 23+\delta})\right\}\frac{1}{(1-c)^{k-1}}\prod_{r=1}^k (1-c)\overset{r-1}{T}. 
\end{split}
\ee 
And consequently, see (\ref{5.2}) and (\ref{5.5}), we obtain the following. 

\begin{mydef45}
\be \label{5.6} 
\int_{T\times \overset{1}{T}\times \dots\times \overset{k-1}{T}}^{[T\times \overset{1}{T}\times \dots\times \overset{k-1}{T}]^1}\left|\zf\right|^2{\rm d}t\sim 
\frac{1}{(1-c)^{k-1}}\prod_{r=1}^k\int_{\overset{r-1}{T}}^{\overset{r}{T}}\left|\zf\right|^2{\rm d}t,\ T\to\infty.  
\ee   
\end{mydef45}

\begin{remark}
We can obtain also the formula: 
\be \label{5.7} 
\int_{(T)^k}^{[(T)^k]^1}\left|\zf\right|^2{\rm d}t\sim\frac{1}{(1-c)^{k-1}}\left\{ \int_T^{\overset{1}{T}}\left|\zf\right|^2{\rm d}t\right\}^k,\ T\to\infty, 
\ee  
where, of course, $\overset{1}{T}=[T]^1$. Since 
\bdis 
T<\overset{1}{T}<\dots<\overset{k-1}{T}
\edis 
the formula (\ref{5.7}) does not follow directly from (\ref{5.6}). 
\end{remark} 

\begin{remark}
Since\footnote{See (\ref{3.4}) with $r=1$.} 
\be \label{5.8} 
\int_{T+\dots+\overset{k-1}{T}}^{[T+\dots+\overset{k-1}{T}]^1}\left|\zf\right|^2{\rm d}t=\sum_{r=1}^{k}(1-c)\overset{r-1}{T}+\mcal{O}(T^{a+\delta}), 
\ee  
then (\ref{1.9}) follows directly. Mixed formulas of type (\ref{1.10}) can be obtained by combining corresponding types (\ref{1.7}) and (\ref{1.9}). 
\end{remark} 

\section{Nonlinear phenomena in which Hardy - Littlewood integral is generated by the increment (\ref{1.4})} 

The increments of the Hardy - Littlewood integral are defined, of course, as differences of integrals, see (\ref{1.4}). Now, we can provide, in opposite direction, a nonlinear formula that generates the Hardy - Littlewood integral by means of the increments. 

\subsection{} 

\begin{mydef12}
The following formula 
\be \label{6.1} 
\begin{split}
& \int_0^{\overset{r-1}{T}}\left|\zf\right|^2{\rm d}t= \\ 
& \frac{1}{1-c}\int_{\overset{r-1}{T}}^{\overset{r}{T}}\left|\zf\right|^2{\rm d}t\times \left\{
\frac{e^{\Lambda_1}}{1-c}\int_{\overset{r-1}{T}}^{\overset{r}{T}}\left|\zf\right|^2{\rm d}t\right\}+ \\ 
& \mcal{O}(T^{a+\delta}\ln T),\ r=1,\dots,k,\ T\to\infty 
\end{split}
\ee 
holds true, where 
\be \label{6.2} 
\Lambda_1=2c-1-\ln 2\pi. 
\ee 
\end{mydef12} 

\begin{remark}
Nonlinear transformation of the increments (for segments $[\overset{r-1}{T},\overset{r}{T}]$) onto the corresponding Hardy-Littlewood integral (on adjacent segments $[0,\overset{r}{T}]$)
\be \label{6.3} 
\int_{\overset{r-1}{T}}^{\overset{r}{T}}\left|\zf\right|^2{\rm d}t\xrightarrow{\mcal{T}_1}\int_0^{\overset{r-1}{T}}\left|\zf\right|^2{\rm d}t,\ r=1,\dots,k
\ee 
is defined by the formula (\ref{6.1}). 
\end{remark} 

\subsection{} 

\begin{mydef13}
The following formula 
\be \label{6.4} 
\begin{split}
& \int_0^{\overset{r}{T}}\left|\zf\right|^2{\rm d}t= \\ 
& \frac{1}{1-c}\int_{\overset{r-1}{T}}^{\overset{r}{T}}\left|\zf\right|^2{\rm d}t\times \left\{
\frac{e^{\Lambda_2}}{1-c}\int_{\overset{r-1}{T}}^{\overset{r}{T}}\left|\zf\right|^2{\rm d}t\right\}+ \\ 
& \mcal{O}(T^{a+\delta}\ln T),\ r=1,\dots,k,\ T\to\infty 
\end{split}
\ee 
holds true, where 
\be \label{6.5} 
\Lambda_2=c-\ln 2\pi. 
\ee 
\end{mydef13} 

\begin{remark}
Nonlinear transformation of the increments (for segments $[\overset{r-1}{T},\overset{r}{T}]$) onto the corresponding Hardy-Littlewood integral (on adjacent segments $[0,\overset{r}{T}]$)
\be \label{6.6} 
\int_{\overset{r-1}{T}}^{\overset{r}{T}}\left|\zf\right|^2{\rm d}t\xrightarrow{\mcal{T}_2}\int_0^{\overset{r}{T}}\left|\zf\right|^2{\rm d}t,\ r=1,\dots,k
\ee 
is defined by the formula (\ref{6.4}). Transformation $\mcal{T}_2$ is then a continuation  of the corresponding increment onto the whole segment $[0,\overset{r}{T}]$. 
\end{remark} 

\section{Next steps in the evolution of the Hardy - Littlewood formula and proofs of Theorems} 

\subsection{} 

Let us remind that we have introduced the Jacob's ladder 
\bdis 
\vp_1(T)=\frac 12\vp(T)
\edis 
in our paper \cite{7} (comp. also \cite{8}), where the function $\vp(T)$ is an arbitrary solution of the noblinear integral equation (also introduced in \cite{7}) 
\bdis 
\int_0^{\mu[x(T)]}\left|\zf\right|^2e^{-\frac{2}{x(T)}t}{\rm d}t=\int_0^T\left|\zf\right|^2{\rm d}t, 
\edis  
where each admissible function $\mu(y)$ generates a solution 
\bdis 
y=\vp(T;\mu)=\vp(T);\ \mu(y)\geq 7y\ln y. 
\edis  
We call the function $\vp_1(T)$ as Jacob's ladder since analogy with the Jacob's dream in Chumash, Bereishis, 28:12. 

\subsection{} 

Next, 83 years after HLI formula (\ref{2.5}), we have prove the following new result, see \cite{7}.  

\begin{mydef92}
The Hardy-Littlewood integral has, in addition to previously known HLI expressions possessing an unbounded error terms at $T\to\infty$, the following infinite number of almost exact representations 
\be \label{7.1}  
\begin{split} 
& \int_0^T\left|\zf\right|^2{\rm d}t=\vp_1(T)\ln\{\vp_1(T)\}+ \\ 
& (c-\ln 2\pi)\vp_1(T)+c_0+\mcal{O}\left(\frac{\ln T}{T}\right), 
\end{split} 
\ee 
where $c_0$ is a constant from the Titschmarsh-Kober-Atkinson formula. 
\end{mydef92} 

\subsection{} 

Further, we have obtained, see \cite{7}, (6.2), the following. 

\begin{mydef93}
\be \label{7.2} 
T-\vp_1(T)\sim (1-c)\pi(T);\ \pi(T)\sim \frac{T}{\ln T},\ T\to\infty, 
\ee  
where the Jacob's ladder can be viewed as the complementary function to the function $(1-c)\pi(T)$ in the sense 
\be \label{7.3} 
\vp_1(T)+(1-c)\pi(T)\sim T,\ T\to\infty. 
\ee 
\end{mydef93} 

\begin{remark}
We have used the following change 
\be \label{7.4} 
\frac{T}{\ln T}\to\pi(T),\ T\to \infty 
\ee 
in the first asymptotic formula (\ref{7.2}), see also \cite{7}, (6.2). This is however completely correct in the asymptotic region by the prime-number law. At the same time, we can give also a motivation for the change (\ref{7.4}). Namely, following\footnote{Comp. \cite{12}, p. 2.} 
\bdis 
\zeta(s)=\exp\left\{ s\int_2^\infty\frac{\pi(x)}{x(x^s-1){\rm d}x}\right\},\ \re(s)>1 
\edis 
we can get that the function $\pi(x)$ is one element of the set of five functions generating the function $\zeta(s)$. 
\end{remark} 

\begin{remark}
Our formula (\ref{7.1}) represents the further step in the evolution of the Hardy - Littlewood integral (1918) for the hypertranscendental Riemann's function $|\zf|^2$. 
\end{remark} 

\subsection{} 

Now, we can make a final step toward the proof of Theorem 1. 

\begin{proof}
By making use of the substitution 
\bdis 
T\to \overset{r-1}{T} 
\edis  
in (\ref{2.5}) and 
\bdis 
T\to\overset{r}{T}
\edis 
in (\ref{7.1})\footnote{See also (\ref{3.2}).} we obtain the following couple of formulas. 
\be \label{7.5} 
\int_0^{\overset{r-1}{T}}\left|\zf\right|^2{\rm d}t=\overset{r-1}{T}\ln \overset{r-1}{T} - (1+\ln 2\pi-2c)\overset{r-1}{T}+\mcal{O}(T^{a+\delta}), 
\ee  
\be \label{7.6} 
\int_0^{\overset{r}{T}}\left|\zf\right|^2{\rm d}t=\overset{r-1}{T}\ln \overset{r-1}{T} - (c-\ln 2\pi)\overset{r-1}{T}+c_0+\mcal{O}\left(\frac{\ln T}{T}\right), 
\ee 
$r=1,\dots,k,\ T\to\infty$, correspondingly. If we subtract (\ref{7.5}) from (\ref{7.6}) we obtain instantly (\ref{3.4}). 
\end{proof} 

\begin{remark}
Our formula (\ref{3.4}) represents a kind of interaction between (\ref{7.5}) and (\ref{7.6}) generated by the Jacob's ladder by means of the property\footnote{See (\ref{3.2}).} 
\bdis 
\overset{r}{T}=\vp_1^{-r}(T),\ T\to\infty. 
\edis 
\end{remark} 

\begin{remark}
In fact, we have used the interaction between the HLI formula and our formula (\ref{7.1}) for the first time in our paper \cite{7}. Namely, the formula (\ref{7.2}) represents the result of this interaction, see \cite{7}, (6.2). This last formula generates all the necessary properties of the set of iterations of the Jacob's ladder, see \cite{9}. 
\end{remark} 

\subsection{} 

Now we give the proof of Theorem 2. 

\begin{proof}
We get from (\ref{7.5}) by (\ref{3.3}) and (\ref{6.2}) that 
\be \label{7.7} 
\begin{split}
& \int_0^{\overset{r-1}{T}}\left|\zf\right|^2{\rm d}t=\overset{r-1}{T}\ln\left\{ e^{\Lambda_1}\overset{r-1}{T}\right\}+\mcal{O}(T^{a+\delta})= \\ 
& \overset{r-1}{T}\ln\left\{ e^{\Lambda_1}\overset{r-1}{T}\right\}\left\{ 1+\mcal{O}\left(\frac{T^{-1+a+\delta}}{\ln T}\right)\right\}. 
\end{split}
\ee  
Next, see (\ref{3.4}), 
\be \label{7.8} 
\begin{split}
& \overset{r-1}{T}=\frac{1}{1-c}\int_{\overset{r-1}{T}}^{\overset{r}{T}}\left|\zf\right|^2{\rm d}t+\mcal{O}(T^{a+\delta})= \\ 
& \frac{1}{1-c}\int_{\overset{r-1}{T}}^{\overset{r}{T}}\left|\zf\right|^2{\rm d}t\times\left\{ 1+\mcal{O}(T^{-1+a+\delta})\right\}, 
\end{split}
\ee  
and by (\ref{7.8}) 
\be \label{7.9} 
\begin{split}
& \ln(e^{\Lambda_1}\overset{r-1}{T})=\ln\left\{ \frac{e^{\Lambda_1}}{1-c}\int_{\overset{r-1}{T}}^{\overset{r}{T}}\left|\zf\right|^2{\rm d}t\right\}+\mcal{O}(T^{-1+a+\delta})= \\ 
& \ln\left\{ \frac{e^{\Lambda_1}}{1-c}\int_{\overset{r-1}{T}}^{\overset{r}{T}}\left|\zf\right|^2{\rm d}t\right\}\times 
\left\{ 1+\mcal{O}\left(\frac{T^{-1+a+\delta}}{\ln T}\right)\right\}. 
\end{split}
\ee 
Now, (\ref{6.1}) follows from (\ref{7.7}) by (\ref{7.8}) and (\ref{7.9}). 
\end{proof} 

\begin{remark}
We can obtain (\ref{6.4}) by a similar way. 
\end{remark}

I would like to thank Michal Demetrian for his moral support of my study of Jacob's ladders.


\begin{thebibliography}{29}
\bibitem{1}  
R. Balasubramanian, An improvement of a theorem of Titschmarsh on mean square of $|\zf|^2$, Proc. Lond. Math. Soc., 3, 540 -- 576 (1978). 
\bibitem{2} 
A. Good, Ein $\Omega$-result f\" ur das quadratische Mittel der Riemannschen Zetafunktion auf der kritischen Linie, Invent. Math. 41 (3), 233 -- 251 (1977). 
\bibitem{3} 
G.H. Hardy, J.E. Littlewood, Contribution to the theory of the Riemann zeta-function and the theory of the distribution of Primes, Acta Math. 41 (1), 119 -- 196, (1918). 
\bibitem{4} 
D.H. Lehmer, On the roots of the Riemann zeta-function, Acta Math. 95, 291 -- 298, (1956).  
\bibitem{5} 
J.E. Littlewood, Researches in the theory of the Riemann $\zeta$-function, Proc. Lond. Math. Soc. (2), 20, (1922). 
\bibitem{6} 
A.E. Ingham, Mean-value theorems in the theory of the Riemann zeta-function, Proc. Lond. Math. Soc. (2), 27, 273 -- 300, (1926). 
\bibitem{7}
J. Moser,
`Jacob's ladders and almost exact asymptotic representation of the Hardy-Littlewood integral`,
Math. Notes 88, (2010), 414-422, arXiv: 0901.3937.
\bibitem{8}
J. Moser,
`Jacob's ladders, the structure of the Hardy-Littlewood integral and some new class of nonlinear integral equations`,
Proc. Steklov Inst. 276 (2011), 208-221, arXiv: 1103.0359.
%
\bibitem{9}
J. Moser, Jacob's ladders, reverse iterations and new infinite set of $L_2$-orthogonal systems generated by the Riemann $\zf$-function, arXiv: 1402.2098v1.  
\bibitem{10} 
J. Moser, Jacob's ladders, interactions between $\zeta$-oscillating systems and $\zeta$-analogue of an elementary trigonometric identity, Proc. Stek. Inst. 299, 189 -- 204, (2017). 
\bibitem{11} 
E.C. Titschmarsh, On van der Corput's method and the zeta-function of Riemann, Quart. J. Math. 5, 195 -- 210, (1934). 
\bibitem{12} 
E.C. Titschmarsh, \emph{The theory of the Riemann zeta-function}, Clarendon Press, Oxford, 1951.   
\end{thebibliography}
\end{document}